\documentclass[10pt]{article}
\usepackage{color}
\usepackage{amssymb}
\usepackage{amsthm,array,amssymb,amscd,amsfonts,latexsym, url}
\usepackage{amsmath}
\usepackage[all]{xy}
\newtheorem{theo}{Theorem}[section]
\newtheorem{prop}[theo]{Proposition}

\newtheorem{lemm}[theo]{Lemma}
\newtheorem{coro}[theo]{Corollary}
\newtheorem{rema}[theo]{Remark}
\newtheorem{Defi}[theo]{Definition}

\title{On  the stable birational type of Fano threefolds}
\author{Claire Voisin}
\date{}

\newfont{\gothic}{eufb10}
\begin{document}

\maketitle
\begin{abstract} We show that,  given a very general cubic threefold $X$, any  smooth cubic threefold stably birational to $X$  must be isomorphic to  $X$. The method applies as well to other   Fano 3-folds, like the  sextic double solid, and the quartic threefold. We also discuss the higher-dimensional case.
 \end{abstract}
\section{Introduction}
Recall that two smooth projective varieties $X$ and $X'$ are stably birational if $X\times \mathbb{P}^r$ is birational to $X'\times  \mathbb{P}^{r'}$, for some integers $r,\,r'$ such that ${\rm dim}\,X+r={\rm dim}\,X'+r'$. A variety $X$ is stably rational if it  is stably birational to a point. The following beautiful result was proved recently by de Gaay Fortman, Engel and Schreieder:
\begin{theo}\label{theodegaay} \cite{SEggay} A very general smooth cubic threefold is not stably rational.
\end{theo}
Note that it had been proved long ago by Clemens and Griffiths \cite{CG} that a smooth cubic threefold is never rational. It has been proved in \cite{tschizh} (using AI) that some smooth cubic threefolds are stably rational. However, it is  proved in \cite{gongxu} (using AI)   that $X\times \mathbb{P}^1$ is  not rational for any smooth cubic threefold $X$. One difficulty in the study of rationality of $X\times \mathbb{P}^r$ with $r\geq 1$, hence of the stable rationality problem, is that the Clemens-Griffiths method can be {\it a priori} applied only in dimension 3, because it uses intermediate Jacobians $J^3(X)$ equipped  with their unimodular polarization, and the latter exists only in dimension $3$. The paper \cite{gongxu} uses the completely new tool invented in \cite{KPPY}, which  is adapted to treat rationality in dimension $4$.

In the paper \cite{shinder}, Shinder initiated the study of a different and interesting problem, namely, given a family of smooth projective algebraic  varieties $\mathcal{X}\rightarrow B$,  whose very general fiber
$\mathcal{X}_b$ is not stably rational, what can be said of the set of pairs of points  $b,\,b'\in B$ such that $\mathcal{X}_{b}$ is stably birational to  $\mathcal{X}_{b'}$?  The main result proved in \cite{shinder} is the following:
\begin{theo}\label{theoshi} \cite{shinder} Assume  $\mathcal{X}\rightarrow B$ is the family of smooth hypersurfaces of degree $d$ in $\mathbb{P}^{n+1}$. Then if one fiber $\mathcal{X}_{b_0}$ is not stably rational, two very general fibers  $\mathcal{X}_{b}$,  $\mathcal{X}_{b'}$ are not  stably birational.
\end{theo}

In other words,  there is a variation of the stable birational type of the fibers, once  one fiber  is not stably rational.  The assumption that  at least one fiber is not stably rational had been established by Schreieder in \cite{schreiederjams} in the range $ d\geq {\rm log}_2n+2,\,n\geq3 $. The very recent paper \cite{hotchkiss}  by Hotchkiss (using AI)  establishes this result for hypersurfaces of degree at least $5$, and similarly the paper \cite{bovnbo} by B\"ohning et al. establishes  this result (using AI) for hypersurfaces of degree at least $4$.  (In the quartic fourfold case, it had been obtained much earlier by Totaro in \cite{totaro}.) In the case of cubic threefolds or quartic threefolds, or quartic or sextic double solids,  this result  is proved  respectively in
\cite{SEggay}, \cite{CTPi}, \cite{voisininvent}, \cite{beauville}. More precisely, these papers establish the non-existence of a  decomposition of the diagonal for the very general considered hypersurface, which implies that it is not stably rational (see Section \ref{secreview} for more explanations). In all these cases,  the assumption of Theorem \ref{theoshi} is satisfied, and   we thus conclude that if $X,\,X'$ are very general smooth cubic threefolds, or  smooth quartic threefolds, or  smooth quartic or sextic double solids, then $X$ is not stably birational to $X'$.

 In the appendix to \cite{shinder}, we gave a different proof of Theorem \ref{theoshi}, which works  assuming  the stronger statement that the very general fiber does not have a Chow   decomposition of the diagonal, and more importantly that this can be  proved by applying the specialization method, exhibiting a nontrivial stable birational invariant on the desingularized special fiber, which disappears on the general fiber, (typically,  a  Brauer class or a unramified cohomology class of positive degree). A similar statement, assuming only that the very general fiber does not have a Chow decomposition of the diagonal, and working in nonzero characteristic, is proved by Schreieder in \cite{schreievar}.

  In this paper, we are going to work under the even stronger assumption that the very general member of the family does not have a {\it cohomological}  decomposition of the diagonal. The non-existence of a cohomological decomposition of the diagonal is proved in \cite{voisininvent} for the very general  quartic double solid,  and, following the same proof and  using \cite{kuznetsova}, one can prove the same statement for the  quartic threefold, as it specializes to a double cover of a quadric threefold ramified along a quartic section.  Following the same method, it is established by  Beauville in \cite{beauville} for the  sextic double solids.  Finally, it is also proved in \cite{SEggay}, by a completely different argument,   for the very general cubic threefold.

Our main result  is  the following  statement.
\begin{theo}\label{theocubic} (a) If $X$ is a very general    cubic threefold,  any smooth cubic threefold $X'$ which is stably birational to $X$ is in fact isomorphic to $X$.

(b) The same statement holds for     quartic threefolds and     sextic  double solids.

(c) If $X$ is a very general    quartic double solid, there are only finitely  many smooth quartic double solids $X'$ which are stably birational to $X$.
\end{theo}
\begin{rema}{\rm If one replaces ``stably birational'' by ``birational'', the conclusion follows from the Clemens-Griffiths method \cite{CG} combined with the generic Torelli theorem which is known for all these Fano threefolds. Indeed, if $X$ is birational to $X'$, then one must have by Clemens and Griffiths an isomorphism of principally polarized abelian varieties
$$J(X)\oplus \bigoplus_iJ(C_i)\cong J(X')\oplus \bigoplus_jJ(C'_j) .$$
 As  $J(X) $ is an  irreducible ppav which is not isomorphic to  the Jacobian of curves (which are also irreducible),  one then gets that
 $J(X)\cong J(X')$ as ppav, and the (generic) Torelli theorem implies that $X\cong X'$. }
\end{rema}

  The proof of Theorem \ref{theocubic} is obtained by exploiting further   the main result of \cite[Appendix]{shinder} (presented here as Theorem \ref{theoappenshi}) in order to prove that, in the situation of Theorem \ref{theocubic}, the intermediate Jacobians of $X$ and $X'$ must be  isogenous.  The conclusion of the proof  relies then   on a stronger  form of the Torelli theorem (see Section \ref{sectorelli}) involving isogeny classes (instead of isomorphism classes) of intermediate Jacobians.

Another application of our arguments compares the stable birational types of different Fano threefolds. For example, we prove the following
\begin{theo}\label{theocompare} (Cf. Corollary \ref{propbirisog}(b)) If $X$ is a very general cubic threefold, quartic threefold, quartic or sextic double solid, and $X'$ is a smooth rationally connected variety stably birational to $X$, then $b_3(X')\geq b_3(X)$.

In particular:

 (a) Let $X$ be  a very general cubic threefold. Then $X$ is not stably birational to any rationally connected projective variety  $X'$ with $b_3(X')<10$. Hence it is not stably birational to any smooth complete intersection Fano variety of dimension $\geq 4$.

(b) Let $X$ be  a very general quartic double solid. Then $X$ is not stably birational to any rationally connected threefold with $b_3(X')<20$, and neither   to any smooth complete intersection Fano variety of dimension $\geq 4$.

\end{theo}

The paper is organized as follows. In Section \ref{secreview}, we discuss the various notions of decompositions of the diagonal, the available criteria and non-existence results. Section \ref{secgenfcat} is devoted to recalling the statement and proof of Theorem \ref{theoappenshi} from \cite[Appendix]{shinder}, expressing stable birational equivalence in terms of decompositions of the diagonal. We establish in that section the following result, from which the other theorems follow directly.

\begin{prop}\label{pronewdu1910intro} (See Proposition \ref{pronewdu1910}.) Let $X$, $X'$ be two rationally connected smooth projective varieties, with ${\rm dim}\,X=3$, which are stably birational. Assume that ${\rm Tors }(H^3(X,\mathbb{Z}))=0$ and that $X$ does not have a cohomological decomposition of the diagonal. Then there exists a  nonzero morphism $\gamma: J^3(X)\rightarrow J^3(X')$ of abelian varieties.
\end{prop}
This morphism $\gamma$ has to be an isogeny when $b_3(X)=b_3(X')$ and $J^3(X) $ is simple, which is the first step in the proof of Theorem \ref{theocubic}.

In Section \ref{sectorelli}, we  state and prove a Torelli theorem for isogenies between intermediate Jacobians of complete intersections Fano threefolds, obtained  as a direct consequence of Donagi's variational Torelli theorem \cite{donagi} in cases (a) and (b) of Theorem \ref{theocubic}. This allows us to deduce, from the existence of a nontrivial isogeny between $J^3(X)$ and $J^3(X')$,   that $X$ is isomorphic to $X'$, concluding the proof of Theorem \ref{theocubic} in cases (a) and (b). The quartic double solid case, for which the Donagi variational Torelli theorem is not known, needs an different argument, which leads to the finiteness statement in case (c).

In Section \ref{secdimfour}, we extend conditionally  our results to dimension $4$ and higher. Although the scheme of proof is exactly the same, starting from dimension 5, we need some extra conditions like the integral Hodge conjecture to be satisfied by the very general member of the family in order to extend the results, see Remark \ref{remafin}. Furthermore, the non-existence of a {\it cohomological} (as opposed to Chow-theoretic) decomposition of the diagonal for the very  general Fano hypersurface does not seem  to have been established in dimension $\geq 4$ in the ranges discussed above. We discuss this last point in Section \ref{secreview}. The higher dimensional results are thus conditional  (see Theorem \ref{theodim4}).

\vspace{0.5cm}

{\bf Thanks.} {\it I thank   the University of Sichuan, where this work  was done, for providing excellent working conditions.}
\section{Decompositions of the diagonal \label{secreview}}
We review in this section various notions of decomposition of the diagonal, and known non-existence results.
\begin{Defi} (i) A smooth projective variety $X$  over an algebraically closed field $K$ has a Chow decomposition of the diagonal if, choosing a point $x\in X(K)$ and  denoting by $\Delta_X\subset X\times X$ the diagonal of $X$,  one can write
$$\Delta_X=X\times x+\Gamma\,\,{\rm in}\,\,{\rm CH}^n(X\times X),$$
where the cycle $\Gamma $ is supported on $D\times X$ for some proper closed algebraic subset $D\subset X$.

(ii) Assuming $K=\mathbb{C}$ (as in the rest of the paper), $X$ has a cohomological  decomposition of the diagonal if, choosing a point $x\in X(\mathbb{C})$ and  denoting by $\Delta_X\subset X\times X$ the diagonal of $X$,  one has an equality of cycle classes
$$[\Delta_X]=[X\times x]+[\Gamma]\,\,{\rm in}\,\,H^{2n}(X\times X,\mathbb{Z}),$$
where the cycle $\Gamma $ is supported on $D\times X$ for some proper closed algebraic subset $D\subset X$.
\end{Defi}

Having a Chow decomposition of the diagonal is equivalent to having trivial universal ${\rm CH}_0$ (see \cite{ACTP}), which means that for any overfield $L\supseteq K$, ${\rm CH}_0(X_L)$ is cyclic, generated by $ x_L$.
As the ${\rm CH}_0$ group is a  stable birational invariant, one gets that stably rational varieties admit a Chow decomposition of the diagonal.
Obviously,  when $K=\mathbb{C}$, having a Chow decomposition of the diagonal implies having a cohomological decomposition of the diagonal, by taking Betti cycles classes.  Ottem exhibited in \cite{ottemexample} a $3$-dimensional rationally connected smooth projective variety having a cohomological decomposition of the diagonal but no Chow decomposition of the diagonal, showing that the converse is not true. This cannot happen in dimension $2$, as proved in \cite{voisinjems}.

Obstructions to the existence of a cohomological decomposition of the diagonal  which are possibly nontrivial for some rationally connected smooth projective complex varieties, are described in \cite{AM}, \cite{CTV},  \cite{voisinjems}, \cite{voisincitrouille}. This is summarized in the following
\begin{prop} \label{propexistcohdec}  If a smooth complex projective variety $X$ has a cohomological decomposition of the diagonal, then

a) The topological Brauer group ${\rm Tors}(H^3(X,\mathbb{Z}))$ vanishes. (This condition is not   satisfied by the desingularization of a certain nodal quartic double solid, as proved in \cite{AM}.)

b)  Integral Hodge classes of degree $4$ on $X$ are algebraic. (This condition is reinterpreted in \cite{CTV} as the vanishing of the degree $3$ unramified cohomology of $X$ with value in $\mathbb{Q}/\mathbb{Z}$. It is proved in \cite{CTOj} to be nontrivial for certain unirational $6$-folds and in \cite{schreiederjams} for some unirational fourfolds.)

c)  $X$ has a universal codimension $2$ cycle parameterized by the intermediate Jacobian $J^3(X)$ (this condition is  not   satisfied by the desingularization of  some nodal quartic double solids, as proved in \cite{voisininvent}).

d) If furthermore ${\rm dim}\,X=3$, the minimal class $\frac{\theta^{g-1}}{(g-1)!}$, $g:=\frac{b_3(X)}{2}$, of the  intermediate Jacobian $J^3(X)$ is algebraic. (This last implication is proved in \cite{voisinjems}. This property  is not  satisfied by the very general cubic threefold, as  proved in \cite{SEggay}.)
\end{prop}

In contrast, there are many more obstructions to the existence of a Chow decomposition of the diagonal.
\begin{prop} \label{propexistchowdec} If a smooth projective variety $X$ over a field $K$  has a Chow  decomposition of the diagonal,
 then

a)  (See \cite{totaro}) One has $H^0(X,\Omega_{X/K}^l)=0$ for any $l>0$. (This condition is automatically satisfied if $X$ is rationally connected in characteristic $0$, but not in nonzero characteristic, see \cite[Proposition 5.10]{kollarbook}.)

 b) (See \cite{CTV}) If $K=\mathbb{C}$,
 all unramified cohomology groups $H^i_{\rm nr}(X, A)$ vanish for $i>0$ and any abelian group $A$.
\end{prop}
\begin{rema}{\rm In b), torsion coefficients $A$ are the interesting ones. Indeed, it is proved in \cite{CTV} that unramified cohomology with $\mathbb{Z}$-coefficients vanishes in positive degree on any smooth projective complex rationally connected variety. }
\end{rema}
\begin{rema}{\rm The first two interesting unramified cohomology groups appearing in b), namely $H^2_{\rm nr} (X,\mathbb{Q}/\mathbb{Z})$ and $H^3_{\rm nr} (X,\mathbb{Q}/\mathbb{Z})$ are in fact (for a rationally connected variety),  the groups appearing respectively in Proposition \ref{propexistcohdec} a) and b). }
\end{rema}

Finally the original form of the  specialization method,  as appearing  in  \cite{voisininvent}, provides the following result.
\begin{theo} \cite{voisininvent} \label{theospec} Let $\mathcal{X}\rightarrow B$ be a flat family of projective complex varieties of dimension $>1$  over a smooth connected  curve. Assume the general fiber $\mathcal{X}_b$ is smooth and one  fiber $\mathcal{X}_0$ has ordinary quadratic  singularities.

a) Assume furthermore that  a (equivalently, any) desingularization $\widetilde{X}_0$ of $\mathcal{X}_0$ does not have a cohomological decomposition of the diagonal. Then the very general fiber $\mathcal{X}_b$ does not have a cohomological decomposition of the diagonal.

b) Assume furthermore that  a (equivalently, any) desingularization $\widetilde{X}_0$ of $\mathcal{X}_0$ does not have a Chow  decomposition of the diagonal. Then the very general fiber $\mathcal{X}_b$ does not have a Chow  decomposition of the diagonal.
\end{theo}
The result we are going to use in this paper is the following
\begin{theo}\label{theonodecompcoh}  Let $X$ be either  a very general cubic or quartic threefold, or a very general quartic or sextic double solid. Then $X$ does not have a cohomological decomposition of the diagonal.
\end{theo}
This statement is the synthesis of results proved in \cite{voisininvent}, \cite{beauville}, \cite{kuznetsova} and \cite{SEggay}. Apart from the cubic threefold case, they are obtained by applying  Theorem \ref{theospec}.  The non-existence of a cohomological decomposition of the diagonal for $\widetilde{X}_0$ can be studied using Proposition \ref{propexistcohdec}. In the case treated in \cite{voisininvent}, the varieties are quartic double solids, the special fiber is the Artin-Mumford double solid, and its desingularization does not admit a cohomological decomposition of the diagonal by Proposition \ref{propexistcohdec} a) because it has a nontrivial topological Brauer group. Thus the very general quartic double solid does not have a cohomological decomposition of the diagonal. The same argument was applied by Beauville  in \cite{beauville} to  the sextic double solid. For the quartic threefold, it suffices, in order to apply  Theorem \ref{theospec}, to find a nodal degeneration $X_0$ of the quartic threefold such that a desingularization $\widetilde{X}_0$ of $X_0$ has a nontrivial topological Brauer class. For this, we first  use the Mori specialization (see \cite[Section V.5]{kollarbook}) of $X$ to a general double cover of a quadric $Q\subset \mathbb{P}^4$  ramified along a quartic section. One can then use \cite{kuznetsova} which precisely constructs  an Artin-Mumford type  nodal degeneration for such double covers.
Finally, the non-existence of a cohomological decomposition of the diagonal for the very general cubic threefold is proved in \cite{SEggay}.

Turning to the  non-existence of a Chow decomposition of the diagonal  for $\widetilde{X}_0$,  this  can be studied using Proposition \ref{propexistchowdec}.
Theorem \ref{theospec} b) has been generalized in \cite{CTPi}, which allows in particular to work with worse singularities on the special fiber and with specialization to nonzero characteristic, as used  in \cite{totaro}, where differential forms on the desingularized  special fiber in nonzero characteristic are used as an obstruction to Chow decomposition of the diagonal on its desingularization, as in Proposition \ref{propexistchowdec} a).

These methods have been also improved by  Schreieder in \cite{schreiederquadricsurfece}, with applications in a much broader numerical range in \cite{schreiederjams}, and also more recently by Hotchkiss (using AI) in  \cite{hotchkiss}, Ottem in \cite{ottem}, and B\"{o}hning-von Bothmer-Pirutka-Zhang in  \cite{bovnbo}, which  proves that very general hypersurfaces of degree $\geq 4$ do not have a Chow decomposition of the diagonal.
However, it is not clear that any of these generalizations allows to disprove the existence of a cohomological decomposition of the diagonal on the general fiber. Schreieder's improved criterion  uses  the existence of  a nontrivial unramified cohomology class   on the desingularization of the special fiber.
  Unfortunately, he does not work with nodal degenerations, so even if the degree of that class is not greater than $3$, it is not clear that one can  deduce from  the non-existence of a cohomological decomposition on the desingularization of the special fiber, the non-existence of a cohomological decomposition on the very general fiber.
  In \cite{schreiederjams}, nonzero unramified cohomology classes of any given  degree are constructed on some unirational conic bundles. In particular, Schreieder constructs unirational fourfolds which do not satisfy the integral Hodge conjecture for degree $4$ Hodge classes. This is an obstruction to a cohomological decomposition of the diagonal. It is not clear however that this can be combined with the specialization method in order to show that the very general quartic or quintic fourfold, or the very general quartic, sextic or octic double covers of $\mathbb{P}^4$ do not admit a cohomological decomposition of the diagonal.

\section{Stably birational varieties \label{secgenfcat}}
Let $X$, $X'$ be smooth projective varieties of the same dimension. We proved in \cite[Appendix]{shinder} the following result, that will be essential for the proof of Theorem \ref{theocubic}.
\begin{theo}\label{theoappenshi}  If $X$ and $X'$ are stably birational of the same dimension $n$, there exist correspondences
$$\Gamma\in {\rm CH}^n(X\times X'),\,\Gamma'  \in {\rm CH}^n(X'\times X)$$
with the following properties
\begin{eqnarray}\label{eqcompos}
\Gamma'\circ \Gamma= \Delta_X+W\,\,{\rm in}\,\,{\rm CH}^n(X\times X),\,\,\Gamma\circ \Gamma'= \Delta_{X'}+W'\,\,{\rm in}\,\,{\rm CH}^n(X'\times X'),
\end{eqnarray}
where the cycle $W$ is supported on $D\times X$ for some proper closed algebraic subset  $D\subset X$, and the cycle $W'$ is supported on $D'\times X'$ for some proper closed algebraic subset  $D'\subset X'$.
\end{theo}
We include the  proof for completeness, and because  the reference  \cite[Appendix]{shinder} is not so easily available.
\begin{proof}[Proof of Theorem  \ref{theoappenshi}]
Let $X$, $X'$ be smooth projective of the same dimension $n$, and let \begin{eqnarray}\label{eqPhiPhiprime}  \Phi: X\times \mathbb{P}^r\dashrightarrow X'\times \mathbb{P}^r
\end{eqnarray}
 be a birational map. Let $\Phi':=\Phi^{-1}:X'\times \mathbb{P}^r\dashrightarrow X\times \mathbb{P}^r$ be its inverse.
Let $$ \Gamma_\Phi\subset X\times  \mathbb{P}^r\times X'\times \mathbb{P}^r,\,\,\Gamma_{\Phi'}\subset X' \times  \mathbb{P}^r \times X\times \mathbb{P}^r$$
be the graphs of $\Phi$, resp. $\Phi'$. We use the same letter to denote the classes of $\Gamma_\Phi$, $\Gamma_{\Phi'}$ in ${\rm CH}^{n+r}(X\times  \mathbb{P}^r \times X'\times \mathbb{P}^r)$.
As $\Phi$, $\Phi'$ are birational and inverse to each other, we have
\begin{eqnarray} \label{eqpremrel} \Gamma_{\Phi'}\circ \Gamma_\Phi=\Delta_{X\times \mathbb{P}^r}+W_0\,\,{\rm in}\,\,{\rm CH}^{n+r}(X\times  \mathbb{P}^r\times X\times  \mathbb{P}^r),
\end{eqnarray}
where the cycle $W_0$ is supported on $D_0\times X\times  \mathbb{P}^r$ for some proper closed algebraic subset $D_0$ of $X\times  \mathbb{P}^r$. (Precisely, we can take for  $D_0$ a divisor such that $\Phi_{\mid (X\times \mathbb{P}^r)\setminus D_0}$ is a morphism and   an isomorphism onto its image in $X'\times \mathbb{P}^r$.)
Let $h$ be the positive  generator of the Chow ring of $\mathbb{P}^r$, and let
$$h_2:={\rm pr}_2^*h\in {\rm CH}^{1}(X \times  \mathbb{P}^r\times X'\times  \mathbb{P}^r),\,h'_4:={\rm pr}_4^*h\in {\rm CH}^{1}(X \times  \mathbb{P}^r\times X'\times \mathbb{P}^r),
$$
where ${\rm pr}_i$, $i=2,\,4$, denotes the  projection of $ X\times  \mathbb{P}^r \times X'\times \mathbb{P}^r $ onto its $i$-th factor. (Here we use the notation $h'_4$ to insist on the fact that it comes from the copy of $\mathbb{P}^r$ that is associated  with $X'$ in (\ref{eqPhiPhiprime}).) We use a similar notation for $X' \times  \mathbb{P}^r\times X\times  \mathbb{P}^r$ where the graph of $\Phi'$ lies.

We now write
\begin{eqnarray} \label{eqdecomp} \Gamma_\Phi=\sum_{i,j}h_2^i{h'_4}^j {\rm pr}_{X\times X'}^*\Gamma_{ij}\,\,{\rm in}\,\,{\rm CH}^{n+r}(X\times \mathbb{P}^r\times X'\times \mathbb{P}^r),
\end{eqnarray}
where $0\leq i\leq r,\,0\leq j\leq r$, and similarly
\begin{eqnarray} \label{eqdecompprime} \Gamma_{\Phi'}=\sum_{i,j}{h'_2}^ih_4^j {\rm pr}_{X'\times X}^*\Gamma'_{ij}\,\,{\rm in}\,\,{\rm CH}^{n+r}(X'\times \mathbb{P}^r\times X\times \mathbb{P}^r).
\end{eqnarray}

\begin{lemm} We have
\begin{eqnarray} \label{eqcompophiphiprime} \Gamma_{\Phi'}\circ \Gamma_\Phi=\sum_{i,k,l}h_2^ih_4^l {\rm pr}_{X\times X}^*(\Gamma'_{kl}\circ \Gamma_{i,r-k})\,\,{\rm in}\,\,{\rm CH}^{n+r}(X\times \mathbb{P}^r\times X\times \mathbb{P}^r).
\end{eqnarray}
\end{lemm}
\begin{proof} To compute  the composition, we have to work  on
$$X\times \mathbb{P}^r\times X'\times \mathbb{P}^r\times X\times \mathbb{P}^r.$$
Following the definition of the composition of correspondences, we need to pull-back the cycles $\Gamma_\Phi$, $\Gamma_{\Phi'}$  from the first four  summands
$X\times \mathbb{P}^r\times X'\times \mathbb{P}^r$ and the last four summands $X'\times \mathbb{P}^r\times X\times \mathbb{P}^r$ respectively, take their intersection, and push-forward the result via the projection onto  $X\times \mathbb{P}^r\times  X\times \mathbb{P}^r$ which forgets the middle term $X'\times \mathbb{P}^r$. Using the developments (\ref{eqdecomp}), (\ref{eqdecompprime}), and applying the same notational convention as above,
one sees that one has to push-forward to $X\times \mathbb{P}^r\times X\times \mathbb{P}^r$  the cycle
$$ \sum_{ijkl} {\rm pr}_{X\times X'}^*\Gamma_{ij}\cdot {\rm pr}_{X'\times X}^*\Gamma'_{kl} h_2^i{h'_4}^{j+k}h_6^l.
$$
on $X\times \mathbb{P}^r\times X'\times \mathbb{P}^r\times X\times \mathbb{P}^r$. It it clear that only the terms with $j+k=r$ will survive, which gives (\ref{eqcompophiphiprime}).
\end{proof}
Coming back to our proof, we next use   the fact that
\begin{eqnarray}\label{eqfacile} {\rm pr}_{X\times X*}(\Delta_{X\times \mathbb{P}^r}\cdot h_2^r)=\Delta_X\,\,{\rm in}\,\,{\rm CH}^n(X\times X),
\end{eqnarray}
where the projections are now those defined on $X\times \mathbb{P}^r\times X\times \mathbb{P}^r$.
Equation (\ref{eqpremrel}), combined with (\ref{eqfacile}), then gives
\begin{eqnarray}\label{eqpoussedu1710} \Delta_X={\rm pr}_{X\times X*}(\Delta_{X\times \mathbb{P}^r}\cdot h_2^r)\\ \nonumber
={\rm pr}_{X\times X*}((\Gamma_{\Phi'}\circ \Gamma_\Phi)\cdot h_2^r)-{\rm pr}_{X\times X*}(W_0\cdot h_2^r)\,\,{\rm in}\,\,{\rm CH}^{n}(X\times   X),
\end{eqnarray}
where we observe that the cycle  $W_1:={\rm pr}_{X\times X*}(W_0\cdot h_2^r)$ is supported on $D_1\times X$ for some proper closed algebraic subset $D_1$ of $X$.

Formula (\ref{eqcompophiphiprime}) now gives us
\begin{eqnarray} \label{eqcompophiphiprimepoussehtrois} {\rm pr}_{X\times X*}((\Gamma_{\Phi'}\circ \Gamma_\Phi)\cdot h_2^r)= \sum_{k} \Gamma'_{kr}\circ \Gamma_{0,r-k}.
\end{eqnarray}
Finally, we observe that ${\rm dim}\,\Gamma_{ij}=n-r+i+j$. It follows that
${\rm dim}\,\Gamma_{0,r-k}=n-k$,
hence that
$${\rm dim}\,\Gamma_{0,r-k}<n\,\,{\rm  if} \,\,k>0.$$ Thus, for $k>0$, the composition $\Gamma'_{kr}\circ \Gamma_{0,r-k}$ is supported on $Z_k\times X$, for some  proper closed algebraic subset $Z_k\subset X$. Hence, (\ref{eqcompophiphiprimepoussehtrois}) can be written as

\begin{eqnarray} \label{eqcompophiphiprimepoussehtroismieux} {\rm pr}_{X\times X*}((\Gamma_{\Phi'}\circ \Gamma_\Phi)\cdot h_2^r)=  \Gamma'_{0r}\circ \Gamma_{0r}+W_2\,\,{\rm in}\,\,{\rm CH}^n(X\times X),\end{eqnarray}
where $W_2$  is supported on $D_2\times X$, for some  proper closed algebraic subset $D_2\subset X$ that contains $Z_k$ for all $k>0$.

Combining (\ref{eqcompophiphiprimepoussehtroismieux}) and (\ref{eqpoussedu1710}) proves the first formula in  (\ref{eqcompos}), with $D=D_1\cup D_2$. The second formula, which will not be used here,  is obtained by exchanging the roles of $X$ and $X'$.
\end{proof}
Letting the correspondences $\Gamma_{\Phi}$, $\Gamma_{\Phi'}$, $W$, $W'$ act on $H^3(X,\mathbb{Z})$ or  $H^3(X',\mathbb{Z})$ (via morphisms of Hodge structures), we get
\begin{coro}\label{corogamma} Under the same assumptions on $X$, $X'$, there exist morphisms
$$\gamma:=\Gamma_{\Phi*}: J^3(X)\rightarrow  J^3(X'),\, w:=W_*: J^3(X)\rightarrow  J^3(X),$$
$$\gamma':= \Gamma_{\Phi'*}: J^3(X')\rightarrow  J^3(X),$$
such that \begin{eqnarray}\label{eqgammaw}\gamma'\circ \gamma ={\rm Id}_{J^3(X)}+w: J^3(X)\rightarrow J^3(X)
\end{eqnarray}
and the endomorphism $w\in {\rm End}(J^3(X))$ is given by the action of a cycle $W\in {\rm CH}^n(X\times X)$ supported on $D\times X$.
\end{coro}

\begin{rema}{\rm The intermediate Jacobian $J^3(X)$ exists as an abelian variety for any smooth projective  rationally connected variety $X$, not only in dimension $3$. Dimension $3$ is used only if we want to introduce a principal polarization on it, but we will not need it.}
\end{rema}

We now prove the main result of this section.
\begin{prop}\label{pronewdu1910} Let $X$ be a smooth projective  rationally connected threefold  over $\mathbb{C}$ and let $X'$ be a smooth projective variety which is stably birational  to $X$. Assume that ${\rm Tors }(H^3(X,\mathbb{Z}))=0$ and that $X$ does not have a cohomological decomposition of the diagonal. Then there exists a  nonzero morphism $\gamma: J^3(X)\rightarrow J^3(X')$ of abelian varieties.
\end{prop}
\begin{proof} Note that, in this statement,  $X'$ does not  have {\it a priori}  the same dimension as $X$, but  we can replace $X$ by $X\times \mathbb{P}^s$ or $X'$ by $X'\times \mathbb{P}^s$ to achieve the equality of dimensions, so that we can apply Theorem \ref{theoappenshi} and Corollary \ref{corogamma}. This  gives us a morphism $\gamma: J^3(X)\rightarrow J^3(X')$ and we just have to prove that $\gamma\not=0$. This is where the assumption that ${\rm dim}\,X=3$ will be used.
If $\gamma=0$, equation (\ref{eqgammaw}) becomes
\begin{eqnarray} \label{eqwhengammazero}  {\rm Id}_{J^3(X)}=-w: J^3(X)\rightarrow J^3(X),
\end{eqnarray}
where the endomorphism $w\in {\rm End}(J^3(X))$ is given by the action of a cycle $W\in{\rm CH}^3(X\times X)$ supported on $D\times X$. It follows from (\ref{eqwhengammazero}) and the fact that $H^{3}(X,\mathbb{Z})$ is torsion-free, that
the cycle $\Delta_X+W\in {\rm CH}^3(X\times X)$  has the property that its cycle class $[\Delta_X]-[W]\in H^6(X\times X,\mathbb{Z})$ has its $(3,3)$-K\"unneth component  equal to zero in ${\rm End}(H^3(X,\mathbb{Z}))$. Using  the fact that $H^1(X,\mathbb{Z})=0$ and  $H^{*}(X,\mathbb{Z})$ is torsion-free, we deduce by K\"unneth decomposition that
\begin{eqnarray}\label{eqwhengammazerodecomp0} [\Delta_X]+[W]\in H^{2*}(X,\mathbb{Z})\otimes  H^{2*}(X,\mathbb{Z}).
\end{eqnarray}
 As $H^{2*}(X,\mathbb{Z})$ is algebraic by \cite{voisininthodge}, we can write (\ref{eqwhengammazerodecomp0}) as
\begin{eqnarray}\label{eqwhengammazerodecomp1} [\Delta_X]+[W]=[Z]\,\,{\rm in}\,\, H^{6}(X\times X,\mathbb{Z}),
\end{eqnarray}
where the cycle $Z\in{\rm CH}^3(X\times X)$ is decomposable, that is, of the form $Z=\sum_i{\rm pr}_1^*Z_i\cdot {\rm pr}_2^*Z'_i$ for some cycles $Z_i,\,Z'_i\in{\rm CH}^*(X)$.
Finally, writing  $Z=\alpha \,X\times {\rm pt}+\sum_{{\rm codim}\,Z_i>0}{\rm pr}_1^*Z_i\cdot {\rm pr}_2^*Z'_i$, formula (\ref{eqwhengammazerodecomp1}) and  the fact that  $W$ is supported on $D\times X$  easily   imply $\alpha=1$, which  gives a cohomological decomposition of the diagonal for $X$,  contradicting our assumption.
\end{proof}
\begin{coro}\label{propbirisog}  (a) Let $X$ be a very general cubic theefold,  and let $X'$ be a cubic threefold which is stably birational to $X$. Then there exists an isogeny $J^3(X)\rightarrow J^3(X')$.

The same statement holds for the very general quartic threefold or quartic or sextic  double solid.

(b) (Cf. Theorem \ref{theocompare}.) Let $X$ be a very general cubic threefold, quartic threefold or quartic or sextic  double solid. If $X'$ is any smooth projective  variety  which is stably birational to $X$, then $b_3(X')\geq b_3(X)$.
\end{coro}
\begin{proof} (a)  By Theorem \ref{theonodecompcoh}, a very general cubic threefold  $X$  does not admit a cohomological decomposition of the diagonal. Furthermore, there is no torsion in $H^3(X,\mathbb{Z})$. Proposition \ref{pronewdu1910} thus provides a nonzero morphism $\gamma:J^3(X)\rightarrow J^3(X')$ of abelian varieties.   As $J^3(X)$ is simple because $X$ is very general, hence has Mumford-Tate group equal to the symplectic group of $H^3(X,\mathbb{Q})$, and ${\rm dim}\,J^3(X')={\rm dim}\,J^3(X)$, we get that $\gamma$ is an isogeny.

The other cases  work in the same way. By Theorem \ref{theonodecompcoh}, none of the  very general  Fano threefolds in the considered families admits  a cohomological decomposition of the diagonal. The other properties used above also hold, because they are  ample hypersurfaces in (weighted) projective space.

We finally prove (b). Assuming $X$ is stably birational to $X'$, we get as above a nonzero  morphism
$$\gamma =\Gamma_*: J^3(X)\rightarrow J^3(X'),$$
which has to be injective  since the intermediate Jacobian of $X$ is as before a simple abelian variety.
It follows that  $b_3(X')\geq b_3(X)$.
\end{proof}

\section{Applications of the  variational  Torelli theorem for hypersurfaces \label{sectorelli}}
We complete in this section the proof of Theorem \ref{theocubic}.
We first explain  the statement and the  proof of the following generic  Torelli theorem   up to isogeny for cubic threefolds, quartic threefolds and sextic double solids, with a weaker version for the quartic double solid. The proof relies in the cubic case  on the Carlson-Griffiths  variational  Torelli theorem for cubic hypersurfaces of dimension $3m$ (see \cite{cagri}) and, in the quartic case (or sextic double solid case), on the (adequately generalized)  Donagi variational Torelli theorem for hypersurfaces (see \cite{coxdonagi}, \cite{donagi}, \cite{dosy}).
For statements (a) and (b) below, the proofs  are due to  Carlson-Griffiths \cite{cagri} and Donagi and collaborators (see \cite{coxdonagi}, \cite{donagi}) except for the observation that these results apply as well to prove a generic Torelli for rational Hodge structures, which is what we need here. The proof in  case  (c) is new.
\begin{theo}\label{theotorelliisog}  (a)  Let $X$ be a very general cubic threefold, and let $X'$ be a smooth cubic threefold. Then, if $J(X)$ is isogenous to $J(X')$, $X$ is isomorphic to $X'$.

(b)  The same statement holds for quartic threefolds, and  sextic double solids.

(c) If $X$  a very general quartic double solid, there are at most finitely many    smooth quartic double solids $X'$ such that $J(X')$ is isogenous to $J^3(X)$.
\end{theo}
\begin{rema} {\rm The similar result for curves of genus $\geq 5$ appears in \cite{bardpiro}, and we follow here  the shorter alternative  argument mentioned   at the end of that paper, and which is due to us.}
\end{rema}
\begin{rema}{\rm The statement above  is weaker for quartic double solids, due to the quartic case in Donagi's list of possible  exceptions to his variational Torelli theorem, see \cite{donagi}.}\end{rema}
Before explaining the proof, we complete the proof of Theorem \ref{theocubic}.
\begin{proof}[Proof of Theorem \ref{theocubic}] Let $X$ be a very general cubic threefold and let $X'$ be a smooth cubic threefold which is stably birational to $X$.  By Proposition \ref{propbirisog}, there exists a nontrivial isogeny $J^3(X)\rightarrow J^3(X')$. By Theorem \ref{theotorelliisog}, one has $X\cong X'$.

Statement b) of  Theorem \ref{theocubic} is proved similarly, using statements b) for both Proposition \ref{propbirisog} and Theorem \ref{theotorelliisog}.

Statement (c) of  Theorem \ref{theocubic} follows similarly from statement (c) of Theorem \ref{theotorelliisog}.
\end{proof}

\begin{proof}[Proof of Theorem \ref{theotorelliisog},(a) and (b)]  If $X$ is a very general cubic threefold, $J^3(X)$ has only one polarization $\theta_X$ (up to a coefficient). It follows that the isogeny
$$\gamma: J^3(X)\rightarrow J^3(X')$$
satisfies $\gamma^*\theta_{X'}=\mu \theta_X$ for some nonzero rational number  $\mu$. Equivalently, the isomorphism of Hodge structures
$$\gamma^*: H^3(X',\mathbb{Q})\cong H^1(J^3(X'),\mathbb{Q})\rightarrow H^1(J^3(X),\mathbb{Q})\cong H^3(X,\mathbb{Q})$$
preserves the intersection pairings, up to a multiplicative coefficient $\mu$.
As already mentioned, the period map for cubic threefolds, which associates to a smooth cubic threefold the isomorphism   class of its intermediate Jacobian equipped with its canonical polarization (or equivalently, its polarized Hodge structure on $H^3(X,\mathbb{Z})$), is an immersion. As $X$ is very general and $J^3(X')$ is isogenous to $J^3(X)$, it follows that $X'$ is also very general and, using the local period map, the isogeny $\gamma$ induces, by considering small deformations of  $X$ and by applying a Baire countability argument,  a holomorphic diffeomorphism
$$\psi: U\cong U'$$
between Euclidean open sets in the moduli space of smooth cubic threefolds, where $0\in U$ parameterizes $X$ and $0'\in U'$ parameterizes $X'$. Furthermore, this diffeomorphism is constructed in such a way that, for any $y\in U$, there is an
isogeny (preserving the polarizations up to the coefficient $\mu$)
$$\gamma_y: J^3(\mathcal{X}_y)\rightarrow J^3(\mathcal{X}'_{\psi(y)}),$$
 where $\mathcal{X}\rightarrow U, \,\,\mathcal{X}'\rightarrow U'$ are the local universal (or Kuranishi) families over $U$ and $U'$ respectively. This family of isogenies translates into
an isomorphism of rational Hodge structures
$$ \gamma_H:  H^3(X_y,\mathbb{Q}) \cong  H^3(X_{\psi(y)},\mathbb{Q})$$
and more generally, an isomorphism of variations of  Hodge structures
which is compatible with the polarizations (up to the coefficient $\mu$).
Differentiating this isomorphism gives an isomorphism of polarized infinitesimal variations of Hodge structures:
\begin{eqnarray}\label{eqdiagram}\xymatrix{& T_{U,0}\ar[d]^{\psi_*}&\rightarrow &{\rm Hom}(H^{2,1}(X), H^{1,2}(X))\ar[d]^{\gamma_H}&
\\
& T_{U',0'}&\rightarrow &{\rm Hom}(H^{2,1}(X'), H^{1,2}(X')),&
}
 \end{eqnarray}
 where the vertical isomorphism $\gamma_H$ on the right is the map induced by $\gamma_H$ via the Hodge decomposition on both spaces  $H^3(X_0,\mathbb{C}) $ and $  H^3(X_{0'},\mathbb{C})$ and the dualities
 $$ H^{1,2}(X)\cong H^{2,1}(X)^*,\,\,H^{1,2}(X')\cong H^{2,1}(X')^* $$
 are given by the polarizations.
 Using the Jacobian description (also due to Carlson and Griffiths \cite{cagri}), it is proved in \cite{cagri} that the data on the top line  determine $X$ (and the data on the bottom line  determine  $X'$ as well), hence we get $X\cong X'$. More precisely,  when $X$ is a smooth cubic $3$-fold, the dual of  the  morphism
$$T_{U,0}\rightarrow {\rm Hom}(H^{2,1}(X), H^{1,2}(X))$$
above (with value in the space of symmetric linear maps),
 identifies, thanks to the Carlson-Griffiths Jacobian formalism \cite{cagri}, to the  map induced by multiplication
$$m_f:S^1\otimes S^1\rightarrow R^2_f,$$
where $S^1=H^0(X,\mathcal{O}_X(1))$, $S^2=H^0(X,\mathcal{O}_X(1))$, $f$ is the defining equation of $X$ and
$R^2_f=S^2/J_f^2$, $J_f^2$ being the Jacobian ideal of $f$ in degree $2$. The kernel of $m_f$ equals, (after symmetrization,) the Jacobian ideal $J_f^2\subset S^2$. The end of the  proof of a) is the Mather-Yau  reconstruction theorem, proved in  \cite{mayau} and improved in \cite{donagi}, saying that the isomorphism class under ${\rm Gl}(S^1)$ of a smooth hypersurface of degree $d$  defined by a homogeneous polynomial equation $f\in S^d:={\rm Sym}^{d}S^1$ is determined by its Jacobian ideal $J_f^{d-1}\subset S^{d-1}$. This proves (a).

 Carlson-Griffiths' theorem above recovering a cubic threefold from its infinitesimal variation of Hodge structure does not apply directly to  the  quartic threefolds, and an important extra ingredient in this case  is provided in \cite{donagi}. If $X$ is a smooth quartic threefold defined by an equation $f=0$, the  morphism
$$T_{U,0}\rightarrow {\rm Hom}(H^{2,1}(X), H^{1,2}(X))$$
given by the infinitesimal variation of Hodge structure of $X$, identifies by Carlson-Griffiths theory to
the multiplication map

\begin{eqnarray} \label{eqmfpour4} m_f:R^4_f\rightarrow {\rm Hom}(R^3_f,R^7_f),\end{eqnarray}
where $S^i:=H^0(\mathbb{P}^4,\mathcal{O}_{\mathbb{P}^4}(i))$, $f$ is the defining equation of $X$, $J^l_f\subset S^l$ is  the Jacobian ideal of $f$ in degree  $l$ and
$R^3_f=S^3/J^3_f$, $R^4_f=S^4/J^4_f$.
The extra  step in Donagi's  method to reconstruct $X$ from the data (\ref{eqmfpour4}) is the symmetrizer lemma (see \cite{dosy} or \cite[Proposition 6.21]{voisinbook}). The symmetrizer lemma indeed allows to reconstruct from (\ref{eqmfpour4}) the multiplication map

\begin{eqnarray} \label{eqmfpour4autredegre} m_f:R^1_f\rightarrow {\rm Hom}(R^3_f,R^4_f).\end{eqnarray}
Precisely, one recovers $R^1_f\subset {\rm Hom}(R^3_f,R^4_f)$ as the set of linear maps $\phi:R^3_f\rightarrow R^4_f$ such that, for any $A,\,B\in R^3_f$,
$$m_f(\phi(A))(B)=m_f(\phi(B))(A)\,\,{\rm in}\,\,R^7_f.$$
Similarly, the symmetrizer lemma will allow to reconstruct successively, in a canonical way:
\begin{eqnarray} \label{eqmfpour4autreautredegre} m_f: R^2_f\rightarrow {\rm Hom}(R^1_f, R^3_f),\,\,m_f: R^1_f\rightarrow {\rm Hom}(R^1_f, R^2_f)
\end{eqnarray}
Having reconstructed canonically (\ref{eqmfpour4autredegre}) and  (\ref{eqmfpour4autreautredegre}), recovering the polynomial structure $S^4\cong{\rm Sym}^4S^1$ and the Jacobian ideal $J^4_f\subset S^4$ works as in the previous case. Thus (b) is proved in the quartic case.

 The similar statement for  sextic double solids is proved  along the same lines, as is mentioned in the introduction of \cite{donagi}, because they can be represented as hypersurfaces in a weighted projective space, so that Carlson-Griffiths' theory applies. Let $X$ be a sextic double solid, ramified over $\mathbb{P}^3$ along a sextic surface defined by an equation $f=0$. Let $R^i_f=S^i/J_f^i$ be the degree $i$ piece of the Jacobian ring  of $f$, where $S^i:=H^0(\mathbb{P}^3,\mathcal{O}_{\mathbb{P}^3}(i))$. The infinitesimal variation of Hodge structure of $X$ is given by the multiplication map
\begin{eqnarray}\label{sexticvhs} R^6_f\rightarrow {\rm Hom}(R^5_f,R^{11}_f).
\end{eqnarray}
   The symmetrizer lemma allows to recover from (\ref{sexticvhs}) the multiplication maps
    \begin{eqnarray}\nonumber S^1\rightarrow {\rm Hom}(R^5_f,R^6_f),\\
  \label{eqmultdeg4} R^4_f=S^4\rightarrow {\rm Hom}(S^1,R^5_f),\\
  \nonumber
   R^3_f=S^3\rightarrow {\rm Hom}(S^1,S^4),\\
   \nonumber
   R^2_f=S^2\rightarrow {\rm Hom}(S^1,S^3),
   \\
   \nonumber
    S^1\rightarrow {\rm Hom}(S^1,S^2).\end{eqnarray}
   It is clear that we can recover in turn from these data the polynomial structures $S^i=R^i_f\cong {\rm Sym}^iS^1$ for $i\leq 4$, and for $i=5$, (\ref{eqmultdeg4}) gives us $J^5_f\subset S^5={\rm Sym}^5S^1$, which  determines $f$ by the Mather-Yau theorem. This proves (b) for the sextic double solid.
   \end{proof}
\begin{proof}[Proof of Theorem \ref{theotorelliisog}, (c)] We follow the same reasoning as above, except that, instead of proving that we can recover the general quartic solid $X$ from the data of  its infinitesimal variation of Hodge structure, we will only prove that there are at most finitely many quartic solids $X'$ with  infinitesimal variation of Hodge structure isomorphic to that of $X$.
Let $X$ be a quartic double solid, ramified over $\mathbb{P}^3$ along a quartic surface defined by an equation $f=0$. Let $S^{\bullet}$ be the homogeneous polynomial ring of $\mathbb{P}^3$ and let $R^i_f=S^i/J_f^i$ be the degree $i$ piece of the Jacobian ring of $f$. The infinitesimal variation of Hodge structure of $X$ is given by the multiplication map
\begin{eqnarray}\label{quarticvhs} R^4_f\rightarrow {\rm Hom}(R^2_f,R^6_f),
\end{eqnarray}
which we will rather represent, using Macaulay duality $R^6_f\cong (R^2_f)^*$ (which corresponds to Serre duality $H^{1,2}(X)\cong H^{2,1}(X)^*$), as
the multiplication map
\begin{eqnarray}\label{quarticvhsdeg2} \mu_f: R^2_f\otimes R^2_f\rightarrow (R^4_f)^*=R^4_f.\end{eqnarray}
Note that $R^2_f=S^2=H^0(\mathbb{P}^3,\mathcal{O}_{\mathbb{P}^3}(2))$.

We will prove the following:
\begin{prop}\label{profinpourquar} Let $f\in S^4$ be generic. Let  $W_f\subset S^4 $ be the set of  elements   $f'\in S^4$ such that there exists an isomorphism of Jacobian rings in degrees $2$ and $4$, i.e. a commutative diagram

\begin{eqnarray}\label{eqdiagmultffprime}\xymatrix{ & R^2_f\otimes R^2_f \ar[d]^{\psi_2\otimes \psi_2}&\stackrel{\mu_f}{\rightarrow}& R^4_f \ar[d]^{\psi_4}& \\
& R^2_{f'}\otimes R^2_{f'}&\stackrel{\mu_{f'}}{\rightarrow}& R^4_{f'}  ,&}
\end{eqnarray}

 for some linear isomorphisms  $\psi_2:R^2_f\rightarrow R^2_{f'}$, $\psi_4:R^4_f\rightarrow R^4_{f'}$. Then the orbit of $f$ under ${\rm Gl}(4)$ is an irreducible component of $W_f$.
\end{prop}

  Proposition \ref{profinpourquar} says that the rational map
  $$S^4/{\rm Gl}(4)\dashrightarrow G_{\rm sym}(19,S^2\otimes S^2)/{\rm Gl}(S^2)$$
   where $G_{\rm sym}(19,S^2\otimes S^2)$ is the Grassmannian of symmetric codimension $19$ subspaces of $S^2\otimes S^2$, which to $f$ associates  the  kernel of the multiplication map $\mu_f$, modulo the action of ${\rm Gl}(S^2) $, has $0$-dimensional general fiber. Combined with the  discussion  above, which says that the datum of this kernel is  exactly the datum of the isomorphism class of the IVHS of $f$, this  implies that there are finitely many $X'$ having IVHS isomorphic to that of $X$. This concludes the proof of Theorem \ref{theotorelliisog}, (c).
\end{proof}
\begin{proof}[Proof of Proposition \ref{profinpourquar}]  We will first exhibit a polynomial $f_0\in S^4$ which defines a singular quartic surface, but satisfies the following properties:

\vspace{0.5cm}

{\it (1) The Jacobian ideal $J^4_{f_0}\subset S^4$ has dimension $4$.

(2) There are finitely many isomorphism classes of  element $g\in S^4$ such that there exists a commutative diagram as in (\ref{eqdiagmultffprime}) with $f=f_0,\,g=f'$, is conjugate to $f_0$ under ${\rm Gl}(4)$.}

\vspace{0.5cm}

We will then discuss whether this is enough to imply Proposition \ref{profinpourquar}.

 The element $f_0$ is constructed as follows. Let $q_1,\,q_2\in S^2$ be two quadratic polynomials. We assume that the two quadrics $Q_i$ defined by $q_i$  intersect transversally along a smooth  elliptic curve $E\subset \mathbb{P}^3$. We also assume that $Q_i$ is singular at a unique point $x_i$, (which  does not belong to $E$, as otherwise $E$ would be singular). We choose $f_0=q_1q_2$. Statement (1) above is the following

\begin{lemm}\label{le12dec} The Jacobian ideal of $f_0$ in degree $4$ has dimension $16$.
\end{lemm}
\begin{proof} This is equivalent to saying that there are no nonzero vector field on $\mathbb{P}^3$ that is tangent to the surface $S_0:=\{f_0=0\}$ along $S_0$. Such a vector field has to be tangent to each component $Q_i:=\{q_i=0\}$, $i=1,\,2$ of $S_0$, hence also to $E$. But there is no nonzero vector field on $\mathbb{P}^3$ which is everywhere tangent to $E$ along $E$.
\end{proof}
Statement (2) above is the following
\begin{lemm}\label{le22dec} There are finitely many isomorphism classes  under ${\rm Gl}(4)$ of elements  $f'\in S^4$  such that there exists a diagram  as in (\ref{eqdiagmultffprime}) with $f=f_0,\,g=f'$.
\end{lemm}
\begin{proof} We observe that the Jacobian ideal of $f_0$ is contained in the ideal sheaf of $E\cup \{x_1,\,x_2\}$ because $f_0$ is singular  along   $E\cup \{x_1,\,x_2\}$. It follows that the quadratic ideal $I_{f_0}\subset {\rm Sym}^2 R^2_{f_0}={\rm Sym}^2S^2$ obtained by symmetrizing
${\rm Ker}\,\mu_{f_0}\subset R^2_{f_0}\otimes R^2_{f_0}$  has for base locus the image of $E\cup \{x_1,\,x_2\}$ under the Veronese embedding of $\mathbb{P}^3$ into $\mathbb{P}^9$. If we have a commutative diagram  as in (\ref{eqdiagmultffprime}) with $f=f_0,\,g=f'$, then ${\rm Ker}\,\mu_{g}$ also has the property that it has a base locus in $\mathbb{P}(R^2_g)$ consisting of the union of an elliptic curve $E'$ isomorphic to $E$ and and two points $y_1,\,y_2$. It follows that the Jacobian ideal of $g$ also  has for base locus (in $\mathbb{P}^3_g$) the union of an elliptic curve $E'$ which becomes projectively isomorphic  to $E$ after the Veronese embedding, and  two points $y_1,\,y_2$. There are only finitely many degree $4$ projectively embedded elliptic curves $E'\subset \mathbb{P}^3$  that become  isomorphic to $E$ after a Veronese embedding of $\mathbb{P}^3$ in $\mathbb{P}^9$. Next, the partial derivatives of $g$ have to vanish along $E'\subset \mathbb{P}^3$, and this easily  implies that $g$ is a quadratic polynomial in the quadratic equations $q'_1,\,q'_2$ defining $E'$. Thus we have $g= q'_1q'_2$, for some quadratic polynomials $q'_i$ vanishing on $E'$. However, the Jacobian ideal of $g$ also vanishes on two extra points $y_i$ which are not on $E'$, which implies that they are singular points of one of the quadrics $Q'_i=\{q'_i=0\}$. They cannot be both on the same quadric $Q'_i$, because otherwise $Q'_i$ would have rank $2$ and $E'$ could not be smooth. So we conclude that both quadrics $Q'_i$ are singular. As there are only finitely many singular quadrics in the pencil $\langle q'_1,\,q'_2\rangle$, we get that there are finitely many possibilities for the isomorphism class of  $g$ under ${\rm Gl}(4)$.
\end{proof}
To finish the proof of Proposition  \ref{profinpourquar}, we finally observe that its statement is equivalent to the fact that, given $f\in S^4$ generic, there are finitely many Veronese embedded $\mathbb{P}^3\cong V_3\subset \mathbb{P}(R^2_f)$ such that the quadratic  ideal $I_{V_3}(2)$ of $V_3$ is  contained in  ${\rm Ker}\,\mu_f$ (or rather its symmetrization in ${\rm Sym}^2 R^2_f$) and the resulting ideal
$$({\rm Ker}\,\mu_f)/I_{V_3}(2)\subset H^0(V_3,\mathcal{O}_{V_3}(2))=H^0(\mathbb{P}^3, \mathcal{O}_{\mathbb{P}^3}(4))$$
is the Jacobian ideal in degree $4$ of some $f'\in =H^0(\mathbb{P}^3, \mathcal{O}_{\mathbb{P}^3}(4))$.
It is then clear that the finiteness result of (2) for $f_0$ implies the rigidity  result of Proposition \ref{profinpourquar} for generic $f$.
\end{proof}

\section{Higher dimension \label{secdimfour}}
Although we  used the intermediate Jacobians in the 3-dimensional case, following the Clemens-Griffiths tradition, most of the arguments (assuming we know that the general hypersurfaces do not have a cohomological decomposition of the diagonal) work as well for fourfolds, working with Hodge structure on primitive cohomology, which is an equivalent datum. We spell-out in this section the case of dimension $4$. In higher  dimension,  the statement will be also conditional to the integral Hodge conjecture, so we  state separately the result in   dimension 4.

\begin{theo} \label{theodim4} (a) Assume that    the very general hypersurface $X$ of degree $d\geq3$ in $\mathbb{P}^5$ does not admit a cohomological decomposition of the diagonal. Then, if $X'$ is any smooth projective rationally connected variety which is stably birational to $X$, one has  $b_{4}(X')\geq b_4(X)$.

(b)  Under the same assumption, if $d=5$, any hypersurface $X'$ of degree $d$ in $\mathbb{P}^5$ that is stably birational to $X$ is isomorphic to $X$.

(c) Statements a) and b) hold for quartic, sextic and octic double covers of $\mathbb{P}^4$.

(d) Under the same assumption as in (a), if $d=3$ or $d=4$,  there at most countably many cubic fourfolds (resp. quartic fourfolds) $X'$ which are stably birational to a given very general cubic fourfold (resp.  quartic  fourfold) $X$.
\end{theo}

\begin{proof} Let $X$ be as above and assume that $X'$ is stably birational to $X$.  We first use Theorem \ref{theoappenshi}, that provides  correspondences
$$\Gamma\in{\rm CH}^4(X\times X'),\,\Gamma'\in {\rm CH}^4(X'\times X)$$ such that
\begin{eqnarray}\label{eqetape1}  \Gamma'\circ \Gamma =\Delta_X+W\,\,{\rm in}\,\,{\rm CH}^4(X\times X),
\end{eqnarray}
where $W$ is supported on $D\times X$ for some proper c:losed algebraic subset $D\subset X$.

We now denote by $\gamma:=[\Gamma]_*: H^4(X,\mathbb{Z})_{\rm prim}\rightarrow H^4(X',\mathbb{Z})$, where $H^4(X,\mathbb{Z})_{\rm prim}$ is the primitive degree $4$ cohomology of $X$, (which agrees by the generalized  Noether-Lefschetz theorem with the transcendental part $H^4(X,\mathbb{Z})_{\rm tr}$ of $H^4(X,\mathbb{Z})$ because $X$ is very general,  and $X$ is not a quadric), and
$\gamma':=[\Gamma']_*: H^4(X',\mathbb{Z})\rightarrow H^4(X,\mathbb{Z})$. Equation (\ref{eqetape1}) gives
\begin{eqnarray} \label{eqetape2}  \gamma'\circ \gamma= {\rm Id}_{H^4(X,\mathbb{Z})_{\rm prim}}+w,
\end{eqnarray}
where $w:=[W]_*:  H^4(X,\mathbb{Z})_{\rm prim}\rightarrow H^4(X',\mathbb{Z})$.

As $X$ is very general, the Hodge structure on $H^4(X,\mathbb{Z})_{\rm prim}$ is simple. It follows that $\gamma=0$ or $\gamma$ is injective.  In the situation of  a), if $\gamma$ is injective,  one concludes that
$$ b_4(X')\geq b_4(X')_{\rm tr}+1\geq b_4(X)_{\rm tr}+1 =b_4(X),$$
which proves a).  In case b), if $\gamma$ is injective, it  has to induce an isomorphism of rational Hodge structures
$$H^4(X,\mathbb{Q})_{\rm prim}\cong  H^4(X',\mathbb{Q})_{\rm prim}.$$
We then conclude by the same arguments  as in the previous section, using the Donagi variational Torelli theorem \cite{donagi}, that $X$ is isomorphic to $X'$. This proves b).

We thus only  have to prove in both cases that $\gamma\not=0$. Assume by contradiction that $\gamma=0$. Then (\ref{eqetape2}) becomes

\begin{eqnarray} \label{eqetape3}  0= {\rm Id}_{H^4(X,\mathbb{Z})_{\rm prim}}+w,
\end{eqnarray}
and $w:=[W]_*:  H^4(X,\mathbb{Z})_{\rm prim}\rightarrow H^4(X',\mathbb{Z})$, where $W\in{\rm CH}^4(X\times X)$ is supported on $D\times X$  for some proper closed algebraic subset $D$ of $X$.
Another way to write (\ref{eqetape3}) is to say that the cycle class $\alpha:=[\Delta_X]+[W]$ has the property that
$$\alpha_{*\mid H^4(X,\mathbb{Z})_{\rm prim}}=0,$$
from which one deduces, using the fact that $H^*(X,\mathbb{Z})$ is torsion free,  that the K\"unneth decomposition   of $\alpha$ writes

\begin{eqnarray} \label{eqetape4} \alpha\in H^0(X,\mathbb{Z})\otimes H^8(X,\mathbb{Z})\oplus H^2(X,\mathbb{Z})\otimes H^6(X,\mathbb{Z})\\ \nonumber
\oplus  \mathbb{Z}\, h^2\otimes h^2\oplus  H^6(X,\mathbb{Z})\otimes H^2(X,\mathbb{Z}\oplus H^8(X,\mathbb{Z})\otimes H^0(X,\mathbb{Z}),
\end{eqnarray}
where $h=c_1(\mathcal{O}_X(1))\in H^2(X,\mathbb{Z})$, so that $(h^2)^{\rm perp}=H^4(X,\mathbb{Z})_{\rm prim}$.
All the cohomology groups $H^{2i}(X,\mathbb{Z}),\,i\not=2$ are genereated by classes of algebraic cycles (the only nontrivial case being the case of $H^6(X,\mathbb{Z})$ which is generated by the class of  a line). Furthermore, the first term in (\ref{eqetape4}) has to be $[X\times {\rm pt}]$ since $\alpha$ acts as ${\rm Id}$ on $[ {\rm pt} \times X]$.
It thus follows from (\ref{eqetape4}) that $[\Delta_X]+[W]-[X\times {\rm pt}]$ is the class of a cycle supported over $D'\times X$, hence that $X$ admits a cohomological decomposition of the  diagonal.

\vspace{0.5cm}

The proof of c) is identical (the Donagi  variational  Torelli theorem holds in all three cases).

For the proof of (d), we use  the infinitesimal Torelli theorem, which guarantees that the local period map is an immersion, and that for a very general cubic (resp.  quartic fourfold) $X$, there are at most countably many cubic fourfolds $X'$ (resp.  quartic fourfolds $X'$)  with an isomorphism
$$ H^4(X',\mathbb{Q})_{\rm prim}\cong H^4(X,\mathbb{Q})_{\rm prim}$$
of rational Hodge structures. Combining this fact with the previous argument proves (d).
\end{proof}
\begin{rema}{\rm In the case of degree $3$ (cubic fourfolds), we know by the structure of the period map that there are countably many cubic fourfolds $X'$ with rational Hodge structure isomorphic to that of a given cubic fourfold $X$, so we cannot improve (d) in that case. In the case  of quartic fourfolds, which is one of the potential exceptions in \cite{donagi}, the variational Torelli theorem we need in order to prove the Torelli theorem with rational coefficients is not known but it is possible, (with a method similar to the one used for  proof of Theorem \ref{theotorelliisog} (c),) to improve statement (d) in that case, proving that there are at most finitely many quartic fourfolds with rational Hodge structure isomorphic to that of a given general quartic fourfold $X$. This also implies by our arguments that, if the very general quartic fourfold has no cohomological  decomposition of the diagonal,  there are finitely many quartic fourfolds stably birational to a very general one.}
\end{rema}
 \begin{rema} \label{remafin} {\rm Starting from dimension $n\geq 5$,  this proof will not work as it is, because we will need the integral cohomology classes $\frac{h^{k}}{d}$ to be algebraic for $k\geq m+1$, where $n=2m$ or $n=2m+1$, and this is not known for general Fano hypersurfaces, that is, in the range  $d\leq n+1$. This statement is known when $X$ is  a cubic of dimension $5$, because it contains a plane, but it is unknown (to our knowledge) for a general quartic or quintic of dimension $5$.}
 \end{rema}

\end{document}